\documentclass[12pt]{amsart}
\usepackage{amsmath,amsfonts,amsthm,amscd,amssymb,euscript,
graphics}%,psfrag}

\newtheorem{theorem}{Theorem}

\newtheorem{question}[theorem]{Question}

\newtheorem{lemma}[theorem]{Lemma}%[section]

\newtheorem{corollary}[theorem]{Corollary}%[section]

\newtheorem{conjecture}[theorem]{Conjecture}

\newtheorem{defn}[theorem]{Definition}%[section]
\newtheorem{proposition}[theorem]{Proposition}%[section]
\newtheorem{varremark}{Remark}
\newenvironment{remark}{\begin{varremark}\em}{\em\end{varremark}}

\newcommand{\Prob}{{\rm Prob}}

\newcommand{\gra}{\mathcal{G}}

\newcommand{\Aut}{\operatorname{Aut}}
\newcommand{\re}[1]{(\ref{#1})}

\newcommand{\beq}{\begin{equation}}
\newcommand{\eeq}{\end{equation}}
\newcommand{\ee}{\end{equation}}

\newcommand{\zed}{{\mathbb Z}}

\newcommand{\reals}{{\mathbb R}}

\newcommand{\fp}{{\mathbb F_p}}

\newcommand{\girth}{{\rm girth}}

\newcommand{\PGL}{{\rm PGL}}
\newcommand{\SL}{{\rm SL}}

\newcommand{\PSL}{{\rm PSL}}
\newcommand{\GL}{{\rm GL}}
\newcommand{\SU}{{\rm SU}}
\newcommand{\SO}{{\rm SO}}

\newcommand{\PSLtwofp}{\PSL_2(\fp)}

\def\mtrx#1#2#3#4{\begin{pmatrix} #1 & #2  \\ #3 & #4\end{pmatrix}}

\newcommand{\SLtwofp}{\SL_2(\fp)}

\newcommand{\Irred} {\mbox{Irred}}

\newcounter{fignum}

\newcommand{\ignore}[1]{}
\newcommand{\eps}{\epsilon}
\newcommand{\Fq}{{\mathbb F}_q}
\newcommand{\Fp}{{\mathbb F}_p}

\title{On the girth of random Cayley graphs}

\author{A. Gamburd}
\address {Department of
Mathematics, University of California, Santa Cruz,  CA 95064 }
\email{agamburd@ucsc.edu}

\author{S. Hoory}
\address{IBM Research Laboratory, Haifa, Israel}
\email{shlomoh@il.ibm.com}

\author{M. Shahshahani}
\address{Institute for Studies in Theoretical Physics and
Mathematics, Tehran, Iran} \email{mehrdads@ipm.ir}

\author{A. Shalev}
\address{Institute of Mathematics, The Hebrew University, Jerusalem, 91904, Israel}
\email{shalev@math.huji.ac.il}

\author{B.  Vir\'ag}
\address{Departments of Mathematics and Statistics, University of Toronto, Ontario, Canada M5S
2E4} \email{balint@math.toronto.edu}

\begin{document}
\bibliographystyle{plain}
\begin{abstract}

We prove that random $d$-regular Cayley graphs of the symmetric
group asymptotically almost surely  have girth at least
$(\log_{d-1}{|G|})^{1/2}/2$ and that random $d$-regular Cayley
graphs of simple algebraic groups over $\Fq$ asymptotically almost
surely  have girth at least $\log_{d-1}{|G|}/\dim(G)$. For the
symmetric $p$-groups the girth is between $\log \log |G|$ and
$(\log|G|)^\alpha$ with $\alpha<1$. Several conjectures and open
questions are presented.
\end{abstract}

\maketitle

\section{Introduction}

The girth of  a graph is the length of a shortest cycle. Finite
regular graphs of large girth are a natural analogue to the
infinite tree. While random regular graphs have nice expansion
properties, their girth tends to be small, as small cycles can
appear at many places independently. The objects of study of this
paper, random Cayley graphs, overcome this problem. While being
random, they are vertex-transitive, giving short cycles fewer
opportunities to appear.

\subsection*{Graphs of large girth}
Let $g = g(n,d)$ be the largest possible girth of a $d$-regular
graph of size at most $n$. Deriving good bounds on $g(n,d)$ for
any $d \geq 3$ is a notoriously hard problem. If we consider $d
\ge 3$ fixed and growing $n$, the best asymptotic estimates known
are:
\begin{eqnarray}\label{gndasymp}
(2+o(1))\cdot \log_{d-1} n \ge g(n,d) \ge (\frac{4}{3}-o(1))\cdot
\log_{d-1} n.
\end{eqnarray}
While it may appear that the problem is essentially solved, the
constant factor gap is crucial here. Clearly, when considering the
inverse of $g$ the constant factor gap becomes an exponent gap.
Also, it is a small miracle that the lower bound constant $4/3$ is
greater than 1, see Conjecture \ref{cayleymoore}.

The first inequality in \re{gndasymp} is a version of the Moore
bound. It is a consequence of a simple counting argument stating
that a ball of radius $\lfloor (g-1)/2 \rfloor$ around a vertex (or
an edge) is a tree, and therefore must have $\Omega((d-1)^{g/2})$
distinct vertices.

For a family of $d$-regular graphs $\gra_i$ of logarithmic girth,
let $\gamma(\{\gra_i\}) =\liminf_{i \to \infty}
\frac{\girth(\gra_i)}{\log_{d-1}(|\gra_i|)}$.  Erd\H os and Sachs
\cite{ES} described a simple procedure yielding families of graphs
with large girth with $\gamma =1$. The first explicit construction
of an infinite degree 4 family with $\gamma
 \approx 0.83$ was given by Margulis \cite{Ma82}, who also gave
 examples of infinite families with arbitrary large degree and
 $\gamma \approx 0.44$; the constructions in question are Cayley
 graphs of $\SLtwofp$. Imrich \cite{im}, extending the work of
Margulis, constructed a family of Cayley graphs of arbitrary
degree with $\gamma \approx 0.48$, and cubic graphs with $\gamma
\approx 0.96$.
 A family of geometrically defined cubic graphs introduced by
 Biggs and Hoare \cite{bh2} was proven to have $\gamma \ge 4/3$ by
 Weiss \cite{weiss}.

 Examples of graphs of arbitrarily large degree
 satisfying $\gamma \ge 4/3$ where given
by Lubotzky,  Phillips and Sarnak \cite{luw}
 and by by Margulis \cite{Ma88}: these are celebrated Ramanujan
 graphs $X^{p,q}$ - Cayley graphs of $PGL_2(q)$ with respect to a
very special  choice of $(p+1)$  generators, where $p$ and $q$ are
primes
 congruent to $-1$ mod $p$ with the Legendre symbol
 $\left(\frac{p}{q}\right)=-1$. A similar result was obtained by Morgenstern \cite{Mo94}
 for any prime power in place of $p$.   Biggs and Bosher \cite{bb} proved that the
 constant $\gamma =\frac{4}{3}$ for $X^{p,q}$ is essentially the
 best possible, namely they showed that
 $$\girth(X^{p, q}) \le 4 \log_{p} q + \log_{p}4 +2.$$
 For every prime power $q$ Lazebnik, Ustimenko, and Woldar \cite{luw}
 constructed families of $q$-regular graphs with $\gamma \ge
 \frac{4}{3} \log_{q}(q-1).$

\subsection*{Random Cayley graphs}
Let $G$ be a finite group and let $S\subset G$. The (undirected)
Cayley graph $\gra(G,S)$ is the undirected graph with the vertex set
$G$ and the edge set $\{(g,gs): g \in G, s \in S\cup S^{-1} \}$.
Given some group $G$, a random $2k$-regular Cayley graph of $G$ is
the Cayley graph $\gra(G,S\cup S^{-1})$ where $S$ is a set of $k$
elements from $G$, selected independently and uniformly at random.

The properties of this model for random graphs received
considerable attention in the last decade. The expansion of such
graphs (for $|S|$ growing with $|G|$) was considered by
Alon-Roichman~\cite{AlRo94}, Pak~\cite{Pa99} and
Landau-Russell~\cite{LaRu04}. The diameter of random Cayley graphs
on the symmetric group was considered by Babai et al.
\cite{BaHe92, BaBeSe04, bh}.

In this work we consider the girth of random Cayley graphs on
various groups. It turns out that random Cayley graphs of the
symmetric group and of the  algebraic groups over finite fields,
tend to have high girth. This is in contrast to random $d$-regular
graphs that tend to have constant girth~\cite{JaLuRu00,McWoWy05}.

\subsection*{Fixed walks in random graphs}

In the classical models for random walks in random environment, an
environment is created by some random process, then a particle
performs a random walk on this environment.

Normally, there are two ways to look at such walks; {\bf quenched}
properties of the random walk are for the typical environment, and
{\bf annealed} properties of the random walk are averaged over
environments.

The model of random Cayley  graphs allows for a third
interpretation, as the random walk (a sequence of symbols
$w=w_1w_2\cdots w_n$ from $S$) can be fixed in advance of generating
the random graphs. Thus, for a random Cayley graph model we can
always talk about the {\bf fixed walk on the random graph}, which,
equivalently, is a random evaluation of the word $w$.

All of our girth results are based on bounds on the {\bf return
probability} of fixed walks on random graphs. More precisely, given
a sequence of groups, we will show that
\begin{equation}\label{return}
\sup_{|w|\le\, \ell_n} P_{G_n} (w=1)=o\left((d-1)^{-\ell_n}\right),
\end{equation}
See Section \ref{words} for further details and discussion.

\subsection*{Results in this paper} We study the girth of random Cayley graphs for three
natural classes of groups. The methods used to prove \re{return} are
unique to each class. While the first-order asymptotics of the girth
is still an open problem in all cases, our results give bounds of
varying precision.

The most general such result is a simple corollary of the
following theorem, due to Dixon, Pyber, Seress and Shalev
\cite{dpss}:

\begin{theorem} \cite{dpss}
Let $G_n$ be a sequence of simple groups with increasing order, and
let $w$ be a word. Then as $n\to \infty$, $P(w=1 \mbox{ in } G_n)\to
0$.
\end{theorem}

\begin{corollary} For $k$ random generators,
$\girth(G_n)\to \infty$ in probability.
\end{corollary}

In the case of the symmetric group, we show the following in
section \ref{sec:sym}.

\begin{theorem}\label{sngirth:theorem}
As $n\to \infty$, a.a.s. the girth of the $d$-regular random Cayley
graph of $S_n$ is at least  $(1/2-o(1)) \cdot \sqrt{\log_{d-1}
|S_n|}$.
\end{theorem}

However, we conjecture that the girth is equal to $O(\log |S_n|)$.

In section \ref{alg} we consider families of simple groups of Lie
type; in this case representations as matrices is very helpful and
we can get stronger bounds.

\begin{theorem} \label{t0}
As $q \to \infty$ a.a.s.  the girth of the $d$-regular random
Cayley graph of  $G(\Fq)$, where $G$ is a simple group of fixed
Lie type and fixed rank over $\Fq$ is at least $(\gamma-o(1))
\log_{d-1} |G(\Fq)|$ with $\gamma = 1/{\dim(G)}$.
\end{theorem}

The bound in Theorem \ref{t0} is optimal except for the crucial
constant $\gamma$. It should be noted that the construction
yielding the lower bound in (\ref{gndasymp}) is a Cayley graph of
$\PGL_2(p)$. However, it seems from computer experiments
(section~\ref{pglcomp:section}), that such a result (or even a
lower bound of $1 \cdot \log_{d-1} |G|$) is unlikely for a random
Cayley graph of $\PGL_2(p)$. In fact, achieving better girth than
$\log_{d-1} |G|$ seems to be a barrier for many combinatorial
constructions such as the result of Erd\H os and Sachs~\cite{ES}.

The $\gamma=1$ threshold can be obtained as a consequence of the
following appealing heuristics. Let $w(x_1, \dots, x_k)$ be a
fixed word in a free group on $k$ generators $x_1, \dots, x_k$.
Let $f_w(g_1, \dots, g_k)$ be an element in $G$ obtained by
substituting $x_i =g_i$.  Define \beq \label{probdef} P_{G}(w) =
\Prob[f_{w}(g_1, \dots, g_k) =1]. \eeq Suppose that for a fixed
short words $w$ we have  $P_{G}(w) \sim 1/|G|$, and that the
events $[f_{w_1} =1]$ and $[f_{w_2}=1]$  are independent for
``generic'' $w_1$ and $w_2$. Then by counting words we could
easily get that girth$(G)/\log_{d-1}|G| \to 1$ as $|G|\to\infty$
a.a.s.

Both of these assumptions are false. It seems that it is not
possible to decrease the satisfaction probability for a sufficiently
large number of words, but it is possible to increase it, for
example in Abelian groups.

The independence assumption is not needed for the lower bound (where
the union bound can be used), but problematic positive correlations
arise when one tries to prove upper bounds. Yet we believe that for
random Cayley graphs, a stronger version of the Moore bound holds
(with constant $1$ instead of $2$).

\begin{conjecture}\label{cayleymoore} Let $\{G_n\}$ be a sequence of groups.
As $n\to \infty$, a.a.s. the girth of the $d$-regular random Cayley
graph of $G_n$ is at most  $(1+o(1))\log_{d-1} |G_n|$.
\end{conjecture}

The third family of groups we are considering are $p$-groups, which
may be thought of as an intermediate class between Abelian groups
(where the girth is at most 4) and simple groups (where the
girth can be logarithmic). This is the only case where we have an
upper bound better than Moore's.

The symmetric $p$-group $W_n(p)$ of height $n$ is the $n$-fold
iterated wreath product of $\zed/(p\,\zed)$. It is
isomorphic to the Sylow $p$-subgroup of the symmetric group
Sym$(p^n)$. It plays the role analogous to the symmetric group in
the realm of $p$-groups: it is a basic family of groups containing
all finite $p$-groups as subgroups.

\begin{theorem}\label{t:symp}
As $n\to \infty$, a.a.s. the girth $g_n$ of the $d$-regular random
Cayley graph of the symmetric $p$-group $G=W_n(p)$ satisfies
$$
(1-o(1))\beta \log\log |G| \le g_n \le (1+o(1))(\log |G|)^{\alpha}
$$
where $\alpha<1$ is a constant depending on $p$ only, and $\beta$
depends on $p$ and $d$.
\end{theorem}

In Section~\ref{sec:genetics}, we present the proof of Theorem
\ref{t:symp}, as well as heuristics to show why the upper bound
should be closer to the truth. It turns out that it helps to
relate this problem to a simple toy model for genetics.

A basic question in this direction is
\begin{question} Does there exist a sequence of $p$-groups of increasing order with
random Cayley graphs of logarithmic girth?
\end{question}

We conclude the paper in section \ref{nondi} by commenting on the
analog of large girth property in the case of compact Lie groups.

\section{The union bound and worst case analysis} \label{words}

In all of our proofs we estimate the probability that a given word
evaluates to the identity, thus creating a short cycle. All of our
girth results are based on bounding the probability of random
elements to satisfy a given word. More precisely, given a sequence
of groups $G_n$, and a word $w$ in $d$ generators for a fixed
$d>2$, we show that
\begin{equation*}
\sup_{|w|\le\, \ell_n} P_{G_n} (w=1)=o\left((d-1)^{-\ell_n}\right),
\end{equation*}
for some sequence $\ell_n$. Summing over all words we get that
\begin{eqnarray}
  P({\rm girth}\, G_n \le \ell_n)&=&
  P(w=1\mbox{ in }G_n\mbox{ for some }|w|\le \ell_n)  \nonumber \\
  &\le & \sum_{|w| \le \ell_n} P(w=1 \mbox{ in }G_n)  \label{union}\\
  &\le&  \# \{w:|w|\le \ell_n\}
  \sup_{|w| \le \ell_n} P(w=1 \mbox{ in }G_n)  \label{sup} \\
  &=& \left(1+d \sum_{l=0}^{l_n}(d-1)^{l}\right)o((d-1)^{\ell_{n}}) \nonumber\\
  &=&o(1) \nonumber
\end{eqnarray}
We believe that the sup bound \re{sup} is wasteful; see Remark
\ref{snpower} in the next section.

The other potentially wasteful part is the union \re{union} bound,
which is not far off when events are not positively correlated.
However, it seems that at least in some cases, there are
correlations. For example, consider words $w$ and $w'$ in two
generators, $a$, $b$. In a significantly large portion of such
words, the exponent sum of $a$ equals 0. Thus $w=w'=1$ if $b=1$,
giving $P(w=w'=1)\ge 1/|G|$. Typically, we expect $P(w=1) \asymp
1/|G|$, and so the uncorrelated case would be $P(w=w'=1)\asymp
1/|G|^2$. We don't know how to take advantage of these correlations
for lower bounds. Moreover, they have blocked our attempts for upper
bounds on girth via the second moment method.

\subsection{Limits of the union bound}

The expression in the union bound \re{union} can be written as a
double sum
$$\sum_{|w| \le \ell_n} P(w=1 \mbox{ in }G_n)  =
|G|^{-k}\sum_{|w|<\ell_n} \sum_{g_1,\ldots, g_k} {\mathbf 1}
(w(g_1,\ldots,g_k)=1)
$$
where $\mathbf{1}$ is 1 if its argument is true and zero otherwise.
Switching the order of summation and changing back to probabilities
gives
$$
|G|^{-k}d(d-1)^{\ell-1}\sum_{g_1,\ldots, g_k} P(w=1)
$$
where $w$ now is a uniform random reduced word of length $\ell$.
If $w$ was just a uniformly chosen word , then $P(w=1)$ would mean
the chance that a random walk on $G$ with generators $\{g_k\}$ is
at the origin at time $n$. It is easy to check (and well-known,
see~\cite{AlSp92} p.139) that for even $\ell$ this probability is
at least $1/|G|$.  Using this it is possible to show that if
$\ell_n \ge (1+\eps) \log_{d-1}|G_n|$, then (\ref{union}) cannot
be $o(1)$, and no proof using the union bound could work to show
that the girth is at least $(1+\eps)\log_{d-1}|G|$.

\subsection{Random evaluation of words}
Bounds on $P(w=1)$ have appeared in the literature \cite{dpss}. A
nice bound, using transitivity properties of groups appears in
\cite{a1}.  In \cite{jones} it is shown that only finitely many
finite simple groups satisfy a given non-trivial law $w$. Word
maps are studied in \cite{larsen, shalev}.

It is also natural to ask (in context of the last paragraph) for
which groups do we have $P(w=1)\ge 1/|G|$. Perhaps surprisingly, It
turns out that this is always true for all words in one or two
generators, but not necessarily for three.  See \cite{a2} for many
counterexamples and discussion.

\section{Random Cayley graphs of $S_n$} \label{sec:sym}

\def\gens{{\sigma_1,\ldots,\sigma_k}}

%In this section we prove Theorem~\ref{sngirth:theorem}.
For some $k \geq 2$, let $\gens$ be independent uniform random
permutations from $S_n$. Let $d=2k$, and $S=\{\sigma_1^{\pm
1},\ldots,\sigma_k^{\pm 1}\}$. We prove for $G=C(S_n,S)$ that
a.a.s $\girth(G) \geq  c \cdot \sqrt{n\log n/\log(d-1)}$ for any
constant $c<1/2$.

\begin{remark}\label{snpower}
The crucial bound used in our proof of
Theorem~\ref{sngirth:theorem} is the upper bound on $P_{G}(w)$,
which holds for all non-trivial words of length at most $l$. We
observe that, for the power word $w=a^l$, this bound is almost
tight. Indeed, $P_{G}(w)$ is at least the probability that the
first $\lfloor n/l \rfloor$ cycles of a random permutation in
$S_n$ have length $l$. Therefore
\[
  P_{G}(w)
  \geq \prod_{i=0}^{\lfloor n/l \rfloor-1} 1/(n-il)
  \geq (1/n)^{\frac n {l}}.
\]
Therefore in order to improve upon Theorem~\ref{sngirth:theorem},
by more than a constant factor, one needs either to avoid the
union bound on $w$ or refine the upper bound on $P_{G}(w)$ to
incorporate more information on the structure of $w$.
\end{remark}

\begin{proof}[Proof of Theorem~\ref{sngirth:theorem}]

Our first observation is that the girth of $G$ is the length of
the shortest non-trivial relation between $\gens$. Therefore,
$\girth(G)\geq g$ with high probability, if for most choices of
$\gens$, no non-trivial word in $\gens$ of length smaller than $g$
is the identity permutation. Clearly, it suffices to check only
non-trivial cyclically reduced words of length $\ell<g$. Namely
words $w = s_0 \cdots s_{\ell-1}$, satisfying $s_i \neq s_{i+1
\pmod \ell}^{-1}$ for $0,\ldots,\ell-1$. We denote the set of such
words by $\Irred_g$; clearly  $|\Irred_g| \leq (d-1)^g$. The
probability of $\girth(G)< g$ is bounded by $\sum_{w \in \Irred_g}
P_{G}(w)$, where $P_{G}(w)$, defined in \eqref{probdef} denotes
the probability that $w$ is the identity permutation.  That is
$P_{G}(w)$ is the probability that $w$ fixes all the $n$ points
$1,\ldots,n$.

Given a word $w = s_0 \cdots s_{\ell-1}$ and some starting point
$x_1$, we trace the path $x_1, x_1 s_0, x_1 s_0 s_1, \ldots$, exposing
the necessary entries of the permutations $\gens$ one by one. In order
that $w$ will fix $x_1$, some coincidence must occur. That is, when
exposing the entries of the path starting at $x_1$, there has to be a
first time when the path arrives at $x_1$ by some permutation
different from $s_0^{-1}$. The probability of such an event occurring
at any specific step $i$ is bounded by $1/(n-e)$, where $e$ is the
number of entries exposed so far.  Since $e$ is at most $\ell$, and
since there are at most $\ell$ choices for $i$, we have $\Pr[x_1 w =
x_1] \leq \ell/(n-\ell)$.

Suppose that we already verified that $w$ fixes
$x_1,\ldots,x_{m-1}$ by exposing the necessary entries. Then we
have exposed at most $(m-1)\cdot \ell$ entries. As long as this
number is smaller than $n$, we can choose a point $x_m$ such that
no entry involving $x_m$ was exposed yet. Repeating the previous
argument, yields an upper bound of $\ell/(n-m\ell)$ on the
probability that $w$ fixes $x_m$, even when conditioning on the
previously exposed entries. Therefore the probability that $w$ is
the identity permutation is bounded by $(\ell/(n-m\ell))^m$, as
long as $m\ell<n$. Substituting $m=n/(2\ell)$, yields the bound
\begin{eqnarray}\label{e:pw-sn}
P(w) \leq (2\ell/ n)^{n/(2\ell)}.
\end{eqnarray}
Therefore
\begin{eqnarray*}
\Pr[\girth(G)<g] \leq |\Irred_g| \cdot (2g/n)^{n/(2g)} \leq
(d-1)^g \cdot (2g/n)^{n/(2g)}.
\end{eqnarray*}
Setting $g=c \cdot \sqrt{n\log n/\log(d-1)}$ for any constant
$c<1/2$, yields the required result:
\begin{eqnarray*}
\Pr[\girth(G)<g]) \leq \exp( -\Omega(\,\sqrt{n\log n\log(d-1)}\,) ).
\end{eqnarray*}

\end{proof}

%%%%%%%%%%%%%%%%%%%%%%%%%%%%%%%%%%%%%%%%%%%%%%%%%%%%%%%%%%%%%%

\section{Random Cayley graphs of simple groups of Lie type} \label{alg}

Before proving Theorem \ref{t0} in general (in section
\ref{pfgen}), we give an elementary proof for the group
$\SL_2(\Fp)$ (in section \ref{sltwoelem}) and discuss computer
experiments (section \ref{pglcomp:section}) and connection between
girth and expansion (section \ref{gexp}) in the case of this
group.

\subsection{Random Cayley graphs of $\SL_2(\Fp)$}
\label{sltwoelem} We begin by giving an elementary proof of the
lower bound on the girth of a random $2k$-regular Cayley graph of
the group $\PGL_2(\fp)$ for prime $p$; the proof for $\SL_2(\fp)$
is similar. The Cayley graph is constructed with respect to the
set $S=\{g_1^{\pm 1},\ldots,g_k^{\pm 1}\}$, where $d=2k$ and $g_1,
\dots, g_k$ are independent uniform random elements from
$\PGL_2(\fp)$.

\begin{theorem}\label{pglgirth:theorem}
As $p\to \infty$, a.a.s.the girth of the $d$-regular random Cayley
graph of $G=\PGL_2(\fp)$ or of $G=\SL_2(\fp)$ is at least
$(1/3-o(1)) \cdot \log_{d-1} |G|$.
\end{theorem}

Before proceeding with the proof of Theorem~\ref{pglgirth:theorem}
we recall the  upper bound on the number of  projective zeros of a
polynomial.

\begin{theorem}[Serre~\cite{Se92}, S{\o}rensen~\cite{So91}]
\label{serre:theorem} Homogeneous polynomial in $m$ variables in
$F_p$ of degree $d$ has at most $dp^{m-2}+(p^{m-2}-1)/(p-1)$
projective zeros, and this is sharp.
\end{theorem}

To prove Theorem~\ref{pglgirth:theorem}, we start with the
following lemma:

\begin{lemma}
Let $w$ be a word of length $\ell$ in the free group $\mathcal
F_k$. If $w$ is not identically 1 for every substitution of values
from $PGL_2(p)$, then for a random substitution
$$
\Pr[w=1]\le \ell/p + O(p^{-2})
$$
where implied constant depends on $k$ only.
\end{lemma}

\begin{proof}
The word $w(g_1,\ldots, g_k)$ evaluated in $\GL_2$ is a matrix
whose entries are rational functions of the entries of the $g_i$.
The reason they are not polynomials is that $w$ may contain
inverses of the form $g_i^{-1}$, so that a factor of $1/\det(g_i)$
appears. Nevertheless, the equation $w=I \times$constant (i.e.
$w=1$ in $\PGL$) reduces to three homogeneous polynomial equations
of degree $\ell$ in $4k$ variables, corresponding to the equations
$a_{11}=a_{22}$, $a_{12}=0$ and $a_{21}=0$.

By our assumptions at least one of these equations is not
identically zero. So by Theorem~\ref{serre:theorem} it has at most
$\ell p^{4k-1}+ O(p^{4k-2})$ solutions among all possible matrices
$g_1,\ldots, g_k$, and therefore there are at most this many in
the subset $\left(\GL_2(p)\right)^k$. Since multiplication by
constant matrices preserves solutions, it follows that there are
at most $\ell p^{4k-1}(p-1)^{-k} + O(p^{3k-2})$ solutions to $w=1$
in $\PGL_2(p)^k$. Dividing by the $k$-th power of
$|\PGL_2(p)|=p(p-1)(p+1)$, completes the proof.
\end{proof}

\begin{proof}[Proof of Theorem~\ref{pglgirth:theorem}]

Let $d=2k$. The number of words of length $\ell$ or less is at
most $(d-1)^{\ell+1}$. The probability of each such word is at
most $\ell/p+O(p^{-2})$. So by the union bound all we  need is
that $(d-1)^{\ell+1}\ell/p=o(1)$, which holds if
\begin{eqnarray*}
\ell=\log_{d-1} p - 2 \log_{d-1} \log_{d-1} p =(1/3-o(1)) \cdot
\log_{2k-1} |PGL_2(p)|.
\end{eqnarray*}

We made the assumption that words of length $\ell$ or less do not
yield the identity for all substitutions; this follows from the
following proposition:
\end{proof}

\begin{proposition} \label{p:short}
For any $k$ the length of the shortest non-trivial word $w(x_1,
\dots, x_k)$ such that $f_{w}(g_1, \dots, g_k) =1$ for all $g_1,
\dots, g_k$ in $\SL(2, \fp)$
%(or over $\PGL(2,q)$)
is at least $\Omega(p/\log p)$.
\end{proposition}

Proposition \ref{p:short} follows from Lemma \ref{rel} and
Corollary \ref{cormtt} proved below.

\begin{lemma} \label{rel}
The length of the shortest non-trivial word $w(x_1, x_2)$ such
that $f_{w}(g_1, g_2) =1$ for all $g_1, g_2$ in $\SL(2, \fp)$ is
at least $p$.
\end{lemma}

\begin{proof}
 Suppose we have a word in two generators $g$,
$h$ and let $g=\mtrx{1}{0}{x}{1}$, $h=\mtrx{1}{x}{0}{1}$. By a
simple inductive argument, for all integers $l_1, k_1, \dots, l_n,
k_n$ we have
$$g^{l_1}h^{k_1}g^{l_2}h^{k_2} \dots g^{l_n}h^{k_n} =
\mtrx{f_{11}}{f_{12}}{f_{21}}{f_{22}+l_1 k_1 \dots l_n k_n
x^{2n}},$$ where $f_{11}, f_{12}, f_{21}, f_{22}$ are polynomials
of degree at most $2n-1$.   If the length of the word is less than
$p$ then all $l_i$ and $k_i$ are less than $p$ in absolute value,
and hence $l_1 k_1 \dots l_n k_n$ is not congruent to zero modulo
$p$. Consequently we have that $(f_{22}+l_1 k_1 \dots l_n k_n
x^{2n})-1$ is a nontrivial polynomial of degree $2n$, which has at
most $2n$ roots.  Since $2n$ is clearly also  less than $p$, there
is choice of $x$ for which the polynomial is not zero modulo $p$;
hence for such $x$ we have that $g^{l_1}h^{k_1}g^{l_2}h^{k_2}
\dots g^{l_n}h^{k_n} \ne \mtrx{1}{0}{0}{1} \mod p$.

Similarly, we obtain that $$g^{l_1}h^{k_1}g^{l_2}h^{k_2} \dots
g^{l_n}h^{k_n}g^{l_{n+1}} = \mtrx{v_{11}}{v_{12}}{v_{21}+l_1 k_1
\dots l_n k_n l_{n+1}x^{2n+1}} {v_{22}},$$ where $v_{11}, v_{12},
v_{21}, v_{22}$ are polynomials of degree at most $2n$; and apply
the preceding argument.
\end{proof}

\begin{lemma}\label{morethantwo:lemma}

Let $\omega$ be a non-empty length $l$ reduced word in the $k$
letters $S = \{ g_1^{\pm 1}, g_2^{\pm 1}, \ldots, g_k^{\pm 1} \}$.
Then for any $k>k' \geq 2$, one can find words
$\omega_1,\ldots,\omega_k$ in the letters $S'= \{ g_1^{\pm 1},
g_2^{\pm 1}, \ldots, g_{k'}^{\pm 1} \}$ so that the word $\omega'$
obtained from $\omega$ by substituting $g_i^{\pm 1}$ by
$\omega_i^{\pm 1}$ for $i=1,\ldots,k$, does not reduce to the
empty word. Moreover, one can find such words with $|\omega_i|
\leq 3+2\log_{2k'-1}/\log l$.
\end{lemma}

\begin{corollary} \label{cormtt}
Let $l$ be the length of the shortest non-trivial word in two
letters over the group $G$. Then the length of the shortest
non-trivial word over $G$ in any number of letters is at least
$\Omega(l/\log l)$.
\end{corollary}

\begin{proof}[Proof of Lemma~\ref{morethantwo:lemma}]
Given the word $\omega$ as above, we set $\omega_i = \omega_{i,L}
x_i \omega_{i,R}$, where $\omega_{i,L}$ and $\omega_{i,R}$ are
uniform independent random reduced words of length $s$, and $x_i$
is chosen from $S'$ so that no cancellations occur in $\omega_i$.
We claim that for a sufficiently large length $s$, the resulting
word $\omega'$ does not reduce to the empty word with probability
greater than zero.

Suppose that $\omega'$ reduces to the empty word. Then, one can
obtain the empty word from $\omega'$ by repeatedly deleting
consecutive pairs of a letter and its inverse. Since, no letter
$x_i$ can cancel until one of the half-words
   $\omega_{i,L}$ or
   $\omega_{i,R}$ cancels, two half words that appear consecutively in the
   expanded
   word must cancel. Since two independent length $s$ reduced words cancel
   with
   probability $2k'(2k'-1)^{-(s-1)}$, and since there are only $l-1$
   consecutive pairs of
   half words in the expanded word $\omega'$, one obtains by union bound
   that $\omega'$
   reduces to the empty word with probability at most
   $l(2k')(2k'-1)^{-(s-1)}$, which is less than one for the claimed
   value of $s$.
\end{proof}

\subsection{Computer experiments  for
$\PGL_2(p)$}\label{pglcomp:section}

\def\mean{\mbox{mean}}
\def\med{\mbox{median}}
\def\std{\mbox{std}}
\def\prob{\mbox{prob}}
\def\ga{\gamma>1}
\def\gb{\gamma>0.8}
\def\gc{\gamma>0.6}
\def\go{\mbox{odd girth}}

In contrast to the permutation group $S_n$ and the iterated wreath
product $W_n$, the size of $\PGL_2(p)$ grows moderately with $p$.
This allows getting some intuition on the asymptotic girth of
$\PGL_2(p)$ from computer experiments. We conducted experiments on
random $4$-regular Cayley graphs over $G=\PGL_2(p)$. The
experiments where conducted for varying primes $p$, where each
experiment was repeated $1000$ times. In each case, our computer
program either returned the length of the shortest cycle, or
announced it to be greater than $30$.

In light of the experimental data, we make the following two
conjectures:
\begin{conjecture} The girth of such graphs is almost surely even.
\end{conjecture}
\begin{conjecture} The girth of such graphs is $(c+o(1))\log_3|G|$,
for some constant $c$, satisfying $1/3<c<1$. Furthermore, a.a.s.
the girth is one of two consecutive even numbers.
\end{conjecture}

We give the following excerpt of our experimental data. As
mentioned, for each value of $p$ we computed the girth of $1000$
random $4$-regular Cayley graphs over $\PGL_2(p)$. In the first
table $n_{\mbox{odd}}$ is the number of graphs with odd girth (out
of thousand).

\begin{eqnarray*}
\begin{array}{l|ccccccccc}
p     &
101 & 331 & 1009 & 2003 & 4001 & 10007 & 20011 & 40009 & 100003\\
\hline n_{\mbox{odd}}& 146 & 138 &   66 &   42 &   22 &    16 & 8
&     7 &      1
\end{array}
\end{eqnarray*}

The second table lists, the number of times each even girth was
attained (out of thousand). The two most abundant values of the
girth where marked in bold. Normalizing these two values by
dividing the girth by $\log_3|\PGL_2(p)|$ yields: $0.85, 0.95$ for
$p=1009$; $0.87, 0.95$ for $p=10007$; and $0.83, 0.89$ for
$p=100003$.

\begin{eqnarray*}
\begin{array}{l|rrrrrrrrrrrr}
\mbox{girth}&\le 10 & 12 &  14 &  16 &  18 &  20 &  22 &  24 &  26 &  28 & 30 & > 30 \\
\hline
p=1009      &    52 & 71 & 111 & {\bf 224} & {\bf 295} & 172 &   9 &   0 &   0 &   0 &  0 &0 \\
p=10007     &     9 &  7 &  18 &  38 &  93 & 198 & {\bf 296} & {\bf 296} &  29 &   0 &  0 &0 \\
p=100003    &     0 &  0 &   1 &   5 &   8 &  39 &  60 & 148 &
{\bf 317} & {\bf 342} & 79 &0
\end{array}
\end{eqnarray*}

\subsection{Girth and expansion} \label{gexp}  In \cite{jbag} it is shown that
if $\Sigma_p$ is a symmetric generating set for
 $\PSLtwofp$ ($p$ prime) such that $\girth(\gra(\PSLtwofp, \Sigma_p)) \ge
 c \log p$, where $c$ is independent of $p$,
 then $\gra(\PSLtwofp, \Sigma_p)$ form a family of expanders.
 Combined with Theorem ~\ref{t0} this implies that
 Cayley graphs of $\PSLtwofp$ are expanders with respect to
generators chosen at random in $\PSLtwofp$.
 The following conjecture, combined with the result in
 \cite{jbag},
 would imply that Cayley graphs
 of $\SLtwofp$ are expanders with respect to \emph{any} choice of
 generators.

 \begin{conjecture}  \label{conj:girth} Suppose
 $\langle \Sigma_p \rangle =\PSLtwofp$.
There is a constant $C$, independent of $p$, satisfying the
following property:  the ball of radius $C$ in the generating set
$\Sigma_p$ contains  two elements $g$, $h$ such that
$\girth(\gra(\PSLtwofp, \{g, h\})) \ge
 \frac{1}{C}\log p$.
 \end{conjecture}

\subsection{Proof of Theorem \ref{t0}} \label{pfgen}
Note first, that, by \cite{LiSh}, almost all  $d$-tuples of
elements in $G(q)$ generate $G$ as $q \to \infty$.

Let $d=2k$. Let $F_k$ be the free group on $x_1, \ldots , x_k$.
For $w \in F_k$ and a group $G$ set
\[
V_w(G) = \{ (g_1, \ldots , g_k): g_i \in G, w(g_1, \ldots , g_k)=1
\}.
\]

 Suppose $G(q)$ is a Chevalley group, coming from the simple
algebraic group $G$.  The set $V_{w}$ is an algebraic set in
$G^k$; by Borel's theorem \cite{borel} for a nontrivial word $w$,
$V_{w}$ is a proper subvariety of $G^k$.

Set $e = \dim{G}$. Then we have
\[
\dim{V_w(G)} \le ke-1.
\]
Note however that $V_w(G)$ may well be a reducible subvariety.

Now, suppose $w$ has length at most $l$. We can view elements of
$G$ as matrices (in a natural way if $G$ is classical, or using a
minimal faithful module if $G$ is exceptional). Then the
requirement $w(g_1, \ldots , g_k) = 1$ translates into polynomial
equations of degree at most $l$ in the matrix entries. Denoting by
$r$ the rank of $G$ the number of such equations is bounded above
by $ar^2$ for some absolute constant $a$. To define $V_w(G)$ over
the affine space of matrices we need to add say $f(r,k)$ fixed
equations defining $G^k$ there.

It is known that an affine  variety $V$ of dimension $D$ defined
by $m$ equations of degrees $\le l$ has at most $l^m (q+1)^D$
$q$-rational points. This follows from Bezout theorem and
intersection theory. Moreover, the same applies if, instead of
taking fixed points of Frobenius, we count solutions to $x^q =
h(x)$, which define the finite twisted groups of Lie type. See
Section 10 of Hrushovski \cite{H} for these facts.

Combining this with the information in the previous paragraph
regarding $V_w$ we obtain
\[
|V_w(G(q))| \le  b_1 l^{ar^2}(q+1)^{ke-1},
\]
where $b_1 =b_1(r,k)$ depends on $r$ and $k$. Since $|G(q)| \sim
q^e$ we have \beq \label{var1}  |V_w(G(q))|/|G(q)^k| \le
bl^{ar^2}/q, \eeq where $b = b(r,k)$ is a constant.

Now noting that the expression on the left-hand side of
\eqref{var1} is the probability of the word $w$ being equal to
identity and applying the union bound of section \ref{words}
completes the proof of Theorem \ref{t0}.

\section{Girth for $p$-groups and toy genetics}\label{sec:genetics}

The symmetric $p$-group $W_n(p)$ of height $n$ is the $n$-fold
iterated wreath product of $\zed/(p\,\zed)$. It is isomorphic to
the Sylow $p$-subgroup of the symmetric group Sym$(p^n)$. Also, it
is isomorphic to the automorphism group of the height $n$ rooted
$p$-ary tree. The group $W_n$ plays the  analogous role to the
symmetric group in the realm of $p$-groups: It is a basic family
of groups containing all finite $p$-groups as subgroups. The size
of this group satisfies $\log_p |W_n(p)| = (p^n-1)/(p-1)$.

In this section, we study the girth of the symmetric $p$-group; we
restrict our attention to $p=2$, as it is conceptually and
notationally more clear. Analogous results hold for other primes
$p$.  The symmetric $2$-group is also the graph automorphism group
of the rooted binary tree of height $n$.
 \newcommand{\crs}[1]{{\mbox{\scriptsize $\,\times_{#1}\,$}}}
Each element $g$ of $W_n=W_n(2)$ can be written as $(g_1,g_2)
{\crs{g}}$, where $g_i \in W_{n-1}$ are elements of the automorphism groups of
the two subtrees $T_1,T_2$ of $T$ with roots at level $1$, and
$\crs{g} \in \zed/(2\,\zed)$ either switches $T_1$ and $T_2$ (active) or equals the
identity (inactive).

We start with a word $w$ in some letters $a,b,\ldots $ and their
inverses $a^{-1}={\tilde a}, b^{-1}={\tilde b}\ldots$. Assume that
the values of $\crs{a},\crs{b},\ldots$ are known. Then $w_1$ and
$w_2$ can be expressed in terms $a_1,a_2,b_1,b_2\ldots $.

It is also clear that if $g$ is a uniform random element in $W_n$,
then $g_1,g_2\in W_{n-1}$ and $\crs{g}$ are independent uniform
choices.

The following toy genetics model describes the way words $w_1$ and
$w_2$ (and their recursive offsprings) are determined.

\subsection*{A toy genetics model} Here we describe a
biologically incorrect model for the genome evolution of a strictly
asexual organism, henceforth referred to as an ``amoeba''.

 The DNA of an
amoeba is a sequence of length $l$ of ``forward'' bases, and their
inverses, or ``backward'' pairs. Backward and forward versions of
the same base cannot be next to each other in the DNA.

At each integer time, each amoeba undergoes fission into two
offspring, and its DNA is inherited as follows. First, two fresh
copies of the DNA are created. Then, ``crossing over'' symbols are
introduced as follows. Each pair of forward and backward bases
introduces its own crossover symbol into the sequence: the forward
alleles after their occurrence; the backwards ones, before.

Each crossover symbol is active or inactive, with equal probability,
independently of others. Crossovers happen at the active symbols.

For example, starting with the word $w=\tilde{a}bcaa{\tilde c}$ (where
$\tilde a=a^{-1}$ denotes the backward pair of $a$) the two fresh
copies are ${\tilde a_1}b_1c_1a_1a_1{\tilde c_1}$ and ${\tilde
a_2}b_2c_2a_2a_2{\tilde c_2}$. With the introduction of the
crossover symbols, the word becomes
\[
  {\tilde a}bcaa{\tilde c} \Rightarrow
  \crs{a}{\tilde a}b\crs{b}c\crs{c}a\crs{a}a\crs{a}\crs{c}{\tilde c}.
\]
Say the random settings activate the symbols $\crs{a}$ and
$\crs{b}$, but not $\crs{c}$. Then the DNA of the two offspring are:
\begin{eqnarray}
 w_1={\tilde a_2}b_2c_1a_1a_2\tilde{c_1},\nonumber\\
 w_2={\tilde a_1}b_1c_2a_2a_1\tilde{c_2}.\label{offspring}
\end{eqnarray}
 We are interested in how fast the DNA diversifies. Call an amoeba
{\bf free} if its DNA consists of all different bases. Starting
from a given DNA $w$, how many generations does it take until a
free amoeba is born? In the above example $l=|\omega|=6$ and
$\omega$ has $3$ different bases. After one generation $\omega_1$
  and $\omega_2$ have $5$ different bases each, so they are not free.

{\bf Heuristic.} Very roughly speaking, in each generation, the
number of bases doubles. Thus within a logarithmic number of steps,
an amoeba should emerge with all different bases in its  DNA.

\begin{conjecture} \label{afterheu}
There exists a constant $c>0$ so that starting with any DNA
configuration of length $n$, the probability that there is an amoeba
at generation $c \log n$ with all different bases in her DNA is at
least 1/2.
\end{conjecture}

In fact, there is a simple conjecture that would imply this and
more. We call an integer-valued function from the space of words a
{\bf complexity function} if it satisfies the following properties.
Note that $w_1$ and $w_2$ denote the random DNA of the offspring as
in \re{offspring}.

\begin{enumerate}
\item $\chi(w_i)\le \chi(w)$ for $i=1,2$
\item $\chi(w)\le 0$ iff $w$ consists of different bases
\item Given $w$ with $\chi(w)\ge 1$ we have
$\Pr[\min(\chi(w_1),\chi(w_2))\le \chi(w)-1]\ge 1/2$.
\end{enumerate}

We define
$$
\bar \chi(\ell)= \sup_{|w|\le \ell} \chi(w)
$$

It is not true that the number of bases in the DNA doubles in each
generation with fixed probability. But we believe that there ``the
log number of bases'' can be replaced by some other function of the
DNA so that we get this behavior.

\begin{conjecture} There exists a complexity function $\chi$ with
$\bar \chi(\ell) \le \beta \log \ell$ for a fixed $\beta\ge 1$ and
all $n\ge 1$.
\end{conjecture}

Here we show that

\begin{lemma}\label{cf} The function $\chi(w)=|w|-m(w)$, where $m(w)$ is the number of distinct
bases used by $w$ is a complexity function. \end{lemma}

Clearly, $\bar\chi(l)=l-1$, and the first two properties are
satisfied. For the third, it suffices to prove the following.

\begin{lemma} The probability that
a given offspring of a given amoeba with some fixed DNA $w$ has at
least one more base in her DNA than her parent (given that the
parent has a repeated base) is at least $1/2$.
\end{lemma}

\begin{proof} We consider a repeated base, say $a$,  for which the two repetitions
are closest to each other in the DNA.

If they have the same orientation, then at the time of fission there
will be a single crossover symbol $\crs{a}$ in between the two.
Given the values of all the other crossover symbols, this symbol is
independent and random, and is active with probability $1/2$. Thus,
in the child, the first occurrences of $a$ are from the same copy of
the DNA or a different copy, with probability $1/2$ each.

If they have different orientation, then there must be at least
another base, say $b$, in between the two occurrences; there, $b$
has to appear a single time, otherwise the pair of $a$'s could not
be closest. Thus a single $b$-crossover symbol $\crs{b}$ appears
between the pair of $a$-s. The proof concludes as in the first case,
except we condition on the value of all crossover symbols but
$\crs{b}$.
\end{proof}

\begin{lemma}\label{compl}
If $\chi$ is a complexity function, then the probability, starting
with DNA $w$, that there is no free amoeba at generation $n$ is at
most $p_1(n,|w|)$, where:
\[
  p_1(n,l)
  = \exp\left(-\frac n 4 \left(1-\frac{2\bar\chi(l)}{n}\right)^2\right).
\]
\end{lemma}

\begin{proof}
Let $\ell=|w|$. We consider the evolution of the DNA $w=w(0)$ let
$w(n+1)$ be the the one of the two children of $w(n)$ with lower
complexity. By property (3) of the complexity function, the process
$\chi(w(n))$ is stochastically dominated by $\bar \chi(\ell)-S_n$,
where $S_n$ is a $\mbox{binomial}(n,1/2)$ random variable, i.e. the sum
of $n$ independent random variables taking the values $\{0,1\}$ with
probability $1/2$ each. For such independent ``coin tosses'', we have
the well-known Chernoff (large deviation) bound
$$P[S_n\le\gamma n/2] \le e^{-n(1-\gamma)^2/4}.$$ Setting $\gamma
n/2=\bar \chi(|w|)$ completes the proof.
\end{proof}

After this brief digression into genetics we turn our attention to
the symmetric 2-groups.

\begin{proposition}
There exists $\beta_k>0$ so that for any word $w$ of length at
most $\ell=\lfloor \beta_k n\rfloor $ in $k$ generators we have
$$P(w=1\mbox{ in }W_n)=o((2k-1)^{\ell}).$$ As a consequence, for $k$
random generators, we have $$\operatorname{girth}(W_n)\ge \beta_k
\log_2 \log_2 |W_n| \qquad a.a.s.$$
\end{proposition}
\begin{proof}
Let $n_0<n$.  By Lemmas \ref{cf} and \ref{compl} with probability
$1-p_1(n_0,\ell)$ there is a free amoeba at generation $n_0$.

If level $n_0$ of the tree is not fixed by $w$, then $w\not=1$ in
$W_{n}$ and we are done. If it is fixed by $w$, then the DNA $w'$ of
the free amoeba describes the action of $w$ on $T'$, one of the
subtrees of height $n-n_0$ rooted at level $n_0$.

Since all bases in $w'$ are different, when the random evaluation of
$w'$ gives a uniform random element of $\Aut(T')$. Thus the
conditional probability of $w=1$ in $W_n$ is at most
\[
  p_2(n-n_0)=P(w'=1\mbox{ in } \Aut(T'))=|\Aut(T_{n-n_0})|^{-1}.
\]
Given some value of $\ell$, we set $n_0,n$ so that both $p_1(n_0,\ell)$
and $p_2(n-n_0)$ are $o((2k-1)^{-\ell})$.  First, we set $n_0 =
\alpha_k \ell$ where $\alpha_k$ is a sufficiently large constant to
make $p_1$ small. It is not difficult to verify that $\alpha_k >
4(\log(2k-1)+1)$ suffices. Second, we take $n-n_0$ sufficiently large
so that $p_2^{-1}=|W_{n-n_0}| \gg (2k-1)^\ell$. Here the situation is
much better, and $n-n_0=\Theta(\log \ell)$ suffices.  Putting the two
bounds together yields the lemma, for any
$\beta_k<[4(\log(2k-1)+1)]^{-1}$.
\end{proof}

It was shown in \cite{av} that typical elements have order
$2^{\alpha n+o(1)}$ with $\alpha<1$. This implies the following.

\begin{proposition}\label{spub}
Even for a single random generator, we have
$\operatorname{girth}(W_n)\le (\log |W_n|)^{\alpha+o(1)}$ a.a.s.
\end{proposition}

So the symmetric two-group gives an interesting example of
intermediate girth groups.  Based on the heuristic argument before
Conjecture~\ref{afterheu}, we have

\begin{conjecture}
For $k$ random generators, there is $\beta=\beta_k$ so that we have
$\operatorname{girth}(W_n) = (\log |W_n|)^{\beta+o(1)}$ a.a.s.
\end{conjecture}
Another conjecture by Ab\' ert and the last author is closely
related to this problem.
\begin{conjecture} Let $w$ be a word of length $2^m$. If $m<n$,
then $P(w=1)<1$ in $W_n$.
\end{conjecture}

The $m<n$ condition is sharp, as $w=a^{2^n}$ is satisfied by all
elements  in $W_n$.

In general, we believe that the power word is the easiest to
satisfy.

\begin{conjecture} Let $w$ be a word of length $2^m$. Then
$P(w=1)\le P(a^{2^m}=1)$  in $W_n$.
\end{conjecture}

If true, this conjecture implies that the length of the shortest
non-trivial word satisfied in $W_n$ is $2^n$.

We leave it for the reader to check that this conjecture implies
that the upper bound in Proposition~\ref{spub} has a lower bound
of the same form (with a different constant $\alpha'$).

\section{Noncommutative diophantine property} \label{nondi} In closing we
mention the continuous analog of the notion of large girth
suitable for elements in the group ring of a compact group. It was
introduced in \cite{GJS99} (with $G=SU(2)$) and called there
noncommutative diophantine property.

\begin{defn}[\cite{GJS99}]\label{def:diophant}
For $k\geq 2$, we say that $g_1,g_2,\ldots,g_k\in G$ satisfy
noncommutative diophantine property if there is a
$D=D(g_1,\ldots,g_k)>0$  such that for any $m\geq 1$ and a word
$W_m$ in $g_1,g_2,\ldots,g_k$ of length $m$ with $W_m \neq e$
(where $e$ denotes the identity in $\SU(2)$) we have \beq
\label{e:5.1} ||W_m - e|| \geq\ D^{-m}. \eeq
\end{defn}
\noindent Here
$$
\left|\left|\left[\begin{array}{cc} a&b\\c&d\end{array}
\right]\right|\right|^2\ =\ |a|^2+|b|^2+|c|^2+|d|^2.
$$

Recall that $\theta \in \reals$ is called diophantine if there are
positive constants $C_1$, $C_2$ such that for all $(k, l)\in
\zed^2$ with $k \ne 0$ we have $|k \theta - l| \geq C_1 k^{-C_2}
$. Equivalently, letting $g = e^{2 \pi \theta} \in \SO(2)$, we may
reexpress  this condition as follows:   $|g^k -1| \geq C'_1
k^{-C'_2}.$ A classical result \cite{kh} asserts that diophantine
numbers  $\theta$ are generic in measure in $\reals$.  Given
diophantine $\theta_1, \dots, \theta_k$ and $g_1=e^{2 \pi
\theta_1}, \dots, g_k=e^{2 \pi \theta_k} \in \SO(2)$,  for any
word $W$ in $g_1, \dots, g_k$ of length $m$ we have  $|W_{m} - 1|
\geq \tilde{C}_1 m^{-\tilde{C}_2}$ for some $\tilde{C}_1,
\tilde{C_2}$. In the case of $\SO(3)$, given  $g_1, \dots, g_k$
generating a free subgroup, a pigeonhole argument shows that for
any $m \ge 1$ there is always a word $W$ in $g_1, g_1^{-1}, \dots,
g_k, g_k^{-1}$ of length at most $m$ satisfying
$$\|W-e\| \le \frac{10}{(2k-1)^{m/6}},$$ so the exponential
behavior in the definition above is the appropriate one.

As was first exploited by Hausdorff, for $G=SU(2)$
\cite{Hausdorff} the relation
$$
W_m(g_1,g_1^{-1},\ldots,g_k,g_k^{-1})\ =\ e
$$
where $W_m$ is a  reduced word of length $m\geq 1$ is not
 satisfied
identically in $G^{(k)}$. Hence the sets
$$
V(W_m)\ :=\ \left\{(g_1,\ldots,g_k)|W_m(g)=e\right\}
$$
are of codimension at least one in $G^{(k)}$.  It follows that
$\cup_{m\geq 1}V(R_m)$ is of zero measure in $G^{(k)}$ and also it
is of the first Baire category in $G^{(k)}$.  Thus the generic
$(g_1,\ldots,g_k)\in G^{(k)}$ (in both senses) generates
 the free group.

This holds quite generally: for  $G$ connected, finite-dimensional
non-solvable Lie group it was proved by D.B.A. Epstein ~\cite{E71}
that for each $k > 0$, and for almost all $k$-tuples $(g_1,
\ldots, g_k)$ of elements of $G$, the group generated by $g_1,
\ldots, g_k$ is free on these $k$ elements.

Now the set of $(g_1,\ldots,g_k)\in G^{(k)}$ for which $\langle
g_1,\ldots,g_k\rangle$ is not free is clearly dense
 in $G^{(k)}$
so it follows easily that the set of
 $(g_1,\ldots,g_k)\in G^{(k)}$ which
are {\em not} diophantine is of the second (Baire) category in
$G^{(k)}$. That is to say the topologically generic
$(g_1,\ldots,g_k)$ is
 free but
not diophantine.   On the other hand in \cite{GJS99} it was proved
that the elements with algebraic number entries are diophantine
and the following conjecture was made:

\begin{conjecture}  \label{c:dio} Generic in the measure sense $(g_1,\ldots,g_k)$ is
diophantine.
\end{conjecture}

 Kaloshin and Rodnianski
\cite{KR00a} established the following result towards conjecture
\ref{c:dio}  for almost every pair $(A,B) \in SO(3) \times SO(3)$
there is a constant $D > 0$ such that for any $n$ and any word
$W_n(A, B)$ of length $n$ in $A$ and $B$  the following weak
diophantine property holds:
$$
\|W_n(A,B)-e \| \geq D^{-n^2}.
$$

\textbf{Acknowledgements}.  A.G. is grateful to Akshay
Venkatesh for interest in this work and insightful remarks.
B.V. thanks Mikl\'os Ab\'ert for many useful discussions.
A.G. was supported in part by DARPA and NSF.  S. H. was
supported in part by a PIMS fellowship.  M.S. was supported
by IPM and was visiting Department of Statistics, Stanford
University, during the preparation of this paper. A.S. was
supported in part by grants from the Israel Science
Foundation and from the Binational Science Foundation
United States-Israel. B.V. was supported by Connaught,
Sloan, and NSERC grants.

\end{document}